\documentclass[12pt]{article}
\usepackage{epsfig,amssymb}
\usepackage{graphics}

\begin{document}
\newcommand{\Gm}{G_\mu}
\newcommand{\Pm}{P_\mu}
\newcommand{\Pim}{\Pi_\mu}
\newcommand{\Hm}{H_\mu}
\newcommand{\Lm}{L_\mu}
\newcommand{\Thm}{\Xi_\mu}
\newcommand{\Xim}{\Xi_\mu}

\newcommand{\ux}{u_{,x}}
\newcommand{\uxx}{u_{,xx}}
\newcommand{\ut}{u_{,t}}
\newcommand{\utt}{u_{,tt}}
\newcommand{\vt}{v_{,t}}
\newcommand{\ppt}{p_{,t}}
\newcommand{\ob}{b}
\newcommand{\cc}{c}
\newcommand{\UU}{U}
\newcommand{\sh}{\sharp}
\newcommand{\fl}{\flat}
\newcommand{\odeg}{\; ; \;}
\newcommand{\bif}{\; , \;}
\newcommand{\bis}{\; ; \;}
\newcommand{\bit}{\; , \;}

\newcommand{\hfillbox}{\hfill $\Box$}
\newcommand{\Gf}{{\Gamma_f}}
\newcommand{\sta}{{\rm sta}}
\newcommand{\ext}{{\rm ext}}
\newcommand{\extu}{{\ext}_{\hspace{-.5cm}{}_{{}_{{}_{u\in\calU}}}}}
\newcommand{\extx}{{\ext}_{\hspace{-.5cm}{}_{{}_{{}_{\x\in\calX}}}}}
\newcommand{\exty}{{\ext}_{\hspace{-.5cm}{}_{{}_{{}_{\y\in\calY}}}}}
\newcommand{\extyd}{{\ext}_{\hspace{-.5cm}{}_{{}_{{}_{\y^*\in\calY*}}}}}
\newcommand{\extud}{{\ext}_{\hspace{-.5cm}{}_{{}_{{}_{u^*\in\calU^*}}}}}
\newcommand{\extua}{{\ext}_{\hspace{-.5cm}{}_{{}_{{}_{u\in\calU_a}}}}}
\newcommand{\exte}{{\ext}_{\hspace{-.5cm}{}_{{}_{{}_{\eps\in\calE}}}}}
\newcommand{\extee}{{\ext}_{\hspace{-.5cm}{}_{{}_{{}_{e\in\calE}}}}}
\newcommand{\exted}{{\ext}_{\hspace{-.5cm}{}_{{}_{{}_{\eps^*\in\calE^*}}}}}
\newcommand{\extxn}{{\ext}_{\hspace{-.5cm}{}_{{}_{{}_{\bx\in\real^n}}}}}
\newcommand{\staua}{{\sta}_{\hspace{-.5cm}{}_{{}_{{}_{u\in\calU_a}}}}}
\newcommand{\stax}{{\sta}_{\hspace{-.5cm}{}_{{}_{{}_{\x\in\calX}}}}}
\newcommand{\staxa}{{\sta}_{\hspace{-.5cm}{}_{{}_{{}_{\x\in\calX_a}}}}}
\newcommand{\stayd}{{\sta}_{\hspace{-.5cm}{}_{{}_{{}_{\y^*\in\calY^*}}}}}
\newcommand{\stayad}{{\sta}_{\hspace{-.5cm}{}_{{}_{{}_{\y^*\in\calY^*_a}}}}}
\newcommand{\bEn}{{\bf E}^{(\eta)}}
\newcommand{\bTn}{{\bf T}^{(\eta)}}

\newcommand{\calTn}{\cal{T}^{(\eta)}}
\newcommand{\calEn}{\cal{E}^{(\eta)}}

\newcommand{\calA}{{\cal{A}}}
\newcommand{\barbA}{\bar{\bf A}}
\newcommand{\calO}{{\cal{O}}}
\newcommand{\calXo}{{\cal{X}}^o}
\newcommand{\calXk}{{\cal{X}}_k}

\newcommand{\U}{U}
\newcommand{\Ub}{\bar{U}}
\newcommand{\Uv}{\check{U}}
\newcommand{\Un}{\hat{U}}
\newcommand{\Ut}{\tilde{U}}
\newcommand{\Vv}{\check{V}}
\newcommand{\Vc}{V^c}
\newcommand{\Vn}{\hat{V}}
\newcommand{\Vb}{\bar{V}}
\newcommand{\Vt}{\tilde{V}}
\newcommand{\Gamv}{\check{\Gamma}}
\newcommand{\Gamn}{\hat{\Gamma}}
\newcommand{\Wb}{\bar{W}}
\newcommand{\Pb}{\bar{P}}
\newcommand{\Fb}{\bar{F}}

\newcommand{\cof}{{\rm cof}}
\newcommand{\adj}{{\rm adj}}
\newcommand{\meas}{{\rm meas}}
\newcommand{ \Ii }{{\partial I}}
\newcommand{ \Iu }{{\partial I_u}}
\newcommand{ \It }{{\partial I_f}}
\newcommand{\barIu}{{\bar{I}_u}}
\newcommand{\barIt}{{\bar{I}_f}}
\newcommand{ \Ji }{J_{\partial I}}
\newcommand{ \Ju }{J_{u}}
\newcommand{ \Js }{J_{\sig}}
\newcommand{ \Jj }{{\cal{J}}}
\newcommand{ \calI }{{\cal{I}}}
\newcommand{ \HH}{{H}}
\newcommand{ \PP}{{P}}
\newcommand{ \LL}{{L}}
\newcommand{ \JJ}{{J}}
\newcommand{ \II}{{\Psi}}
\newcommand{\calLs}{{L_{\mbox{\sta}}}}

\newcommand{ \id }{{i_d}}
\newcommand{\bell}{{\mbox{\large $\ell$}}}
\newcommand{\vsig}{\varsigma}
\newcommand{\vrho}{\varrho}
\newcommand{\vsigbar}{\bar{\varsigma}}
\newcommand{\barvsig}{\bar{\varsigma}}
\newcommand{\vsigb}{\bar{\varsigma}}

\newcommand{\pp}{p}
\newcommand{\vv}{v}
\newcommand{\Ww}{\bar{W}}
\newcommand{\ww}{w}
\newcommand{\barww}{\bar{w}}
\newcommand{\Ll}{{\Large \ell}}
\newcommand{\barbw}{\bar{\bf w}}
\newcommand{\uu}{u}
\newcommand{\zz}{z}
\newcommand{\barpp}{\bar{p}}
\newcommand{\barss}{\bar{s}}
\newcommand{\barx}{\bar{x}}
\newcommand{\bary}{\bar{y}}
\newcommand{\qq}{q}
\newcommand{\barqq}{\bar{\qq}}

\newcommand{\jj}{\mbox{\large{$\jmath$}}}
\newcommand{\cl}{\mbox{cl}}

\newcommand{\Lamo}{{\Lam^o}}
\newcommand{\Bo}{{B^o}}
\newcommand{\pol}{\; ; \; }
\newcommand{\oLam}{\Lamo^*}
\newcommand{\calUo}{{\calU^o}}
\newcommand{\Ao}{{A^o}}
\newcommand{\uo}{{u^o}}
\newcommand{\baruo}{{\bar{u}^o}}
\newcommand{\barou}{{\bar{u}^{o*}}}

\newcommand{\Co}{{C^o}}
\newcommand{\ocalX}{\calX^{o*}}
\newcommand{\ocalU}{\calUo^*}
\newcommand{\ou}{u^{o*}}
\newcommand{\infc}{\stackrel{+}{\vee}}
\newcommand{\rarw}{\rightarrow}
\newcommand{\rarwr}{\rightarrow \real}
\newcommand{\rawu}{\rightharpoonup}
\newcommand{\rawuast}{\stackrel{\ast}{\rightharpoonup}}
\newcommand{\pl}{{\parallel}}
\newcommand{\del}{{\delta}}
\newcommand{\barpartial}{{\bar{\partial}}}
\newcommand{\real}{{\mathbb R}} 
\newcommand{\bolM}{{\mathbb M}} 
\newcommand{\boI}{{\mathbb I}}
\newcommand{\boJ}{{\mathbb J}}

\newcommand{\comp}{{\mathbb C}}
\newcommand{\field}{{\mathbb F}}
\newcommand{\bareal}{\bar{\real}}
\newcommand{\realr}{\stackrel{\rightharpoonup}{\real}}
\newcommand{\reall}{\stackrel{\leftharpoonup}{\real}}
\newcommand{\rareal}{{\rightarrow \real}}
\newcommand{\rabareal}{{\rightarrow \bareal}}
\newcommand{\kap}{{\kappa}}
\newcommand{\dotp}{\dot{p}}
\newcommand{\dotkap}{\dot{\kappa}}
\newcommand{\dotkappa}{\dot{\kappa}}
\newcommand{ \doteps}{\dot{\epsilon}}
\newcommand{ \dottau}{\dot{\tau}}
\newcommand{ \doteta}{\dot{\eta}}
\newcommand{ \dotzet}{\dot{\zeta}}
\newcommand{\dotepsilon}{\dot{\epsilon}}
\newcommand{\half}{\frac{1}{2}}
\newcommand{\Ss}{{\bf \Sigma }}
\newcommand{\biota}{\mbox{\boldmath$\iota$}}
\newcommand{\syst}{{\mathbb S}} 
\newcommand{\systo}{{\syst^o}}
\newcommand{\bsig}{\mbox{\boldmath$\sigma$}}
\newcommand{\bvsig}{\mbox{\boldmath$\varsigma$}}
\newcommand{\bmu}{\mbox{\boldmath$\mu$}}
\newcommand{\barbmu}{\bar{\bmu}}
\newcommand{\brho}{\mbox{\boldmath$\rho$}}
\newcommand{\veps}{\xi}
\newcommand{\barveps}{\bar{\xi}}
\newcommand{\bveps}{\mbox{\boldmath$\xi$}}
\newcommand{\bfeta}{\mbox{\boldmath$\eta$}}
\newcommand{\barbsig}{\bar{\mbox{\boldmath$\sigma$}}}
\newcommand{\barbvsig}{\bar{\mbox{\boldmath$\varsigma$}}}
\newcommand{\barbpsi}{\bar{\mbox{\boldmath$\psi$}}}
\newcommand{\beps}{\mbox{\boldmath$\epsilon$}}
\newcommand{\barbeps}{\bar{\mbox{\boldmath$\epsilon$}}}

\newcommand{\bPhi}{{\mbox{\boldmath$\Phi$}}}
\newcommand{\bPsi}{{\mbox{\boldmath$\Psi$}}}
\newcommand{\btau}{{\mbox{\boldmath$\tau$}}}
\newcommand{\bbta}{{\mbox{\boldmath$\beta$}}}
\newcommand{\blam}{{\mbox{\boldmath$\lambda$}}}
\newcommand{\bLam}{{\mbox{\boldmath$\Lambda$}}}
\newcommand{\barbtau}{\bar{\mbox{\boldmath$\tau$}}}
\newcommand{\bchi}{{\mbox{\boldmath$\chi$}}}
\newcommand{\bchio}{\bchi^o}
\newcommand{\barbchio}{\bar{\bchi}^o}
\newcommand{\barbchi}{\bar{\mbox{\boldmath$\chi$}}}
\newcommand{\bxi}{{\mbox{\boldmath$\xi$}}}
\newcommand{\bphi}{{\mbox{\boldmath$\phi$}}}
\newcommand{\barphi}{\bar{\phi}}
\newcommand{\barbphi}{\bar{\mbox{\boldmath$\phi$}}}
\newcommand{\bpsi}{{\mbox{\boldmath$\psi$}}}
\newcommand{\barbxi}{\bar{\mbox{\boldmath$\xi$}}}
\newcommand{\bnu}{{\mbox{\boldmath$\nu$}}}
\newcommand{\barI}{{\bar{I}}}
\newcommand{\barK}{{\bar{K}}}
\newcommand{\barU}{{\bar{U}}}
\newcommand{\barV}{{\bar{V}}}
\newcommand{\hatV}{{\hat{V}}}
\newcommand{\barW}{{\bar{W}}}
\newcommand{\hatW}{{\hat{W}}}
\newcommand{\hatau}{{\hat{\tau}}}
\newcommand{\hatsig}{{\hat{\sigma}}}
\newcommand{\barv}{{\bar{v}}}
\newcommand{\baru}{\bar{u}}
\newcommand{\barh}{\bar{h}}
\newcommand{\bars}{\bar{s}}
\newcommand{\barp}{\bar{p}}
\newcommand{\barq}{\bar{q}}
\newcommand{\la}{\langle}
\newcommand{\ra}{\rangle}
\newcommand{\barho}{\bar{\rho}}
\newcommand{\bartheta}{\bar{\theta}}
\newcommand{\barbe}{\bar{\bf e}}
\newcommand{\barbeta}{\bar{\beta}}
\newcommand{\bareps}{\bar{\epsilon}}
\newcommand{\barmu}{\bar{\mu}}

\newcommand{\barsig}{\bar{\sigma}}
\newcommand{\barchi}{\bar{\chi}}
\newcommand{\barM}{\bar{M}}
\newcommand{\barxi}{\bar{\xi}}
\newcommand{\bartau}{\bar{\tau}}
\newcommand{\barpsi}{\bar{\psi}}
\newcommand{\ba}{{{\bf a}}}
\newcommand{\br}{{{\bf r}}}
\newcommand{\bp}{{{\bf p}}}
\newcommand{\Oo}{\Omega}
\newcommand{\Oot}{{\Omega_t}}
\newcommand{\Gxi}{{\Gamma_\chi}}
\newcommand{\Gt}{{\Gamma_t}}
\newcommand{\Gu}{{\Gamma_u}}
\newcommand{\bE}{{\bf E}}
\newcommand{\bQ}{{\bf Q}}
\newcommand{\ep}{{\epsilon}}
\newcommand{\bM}{{\bf M}}
\newcommand{\bK}{{\bf K}}
\newcommand{\eba}{\begin{array}}
\newcommand{\eea}{\end{array}}
\newcommand{\ebe}{\begin{eqnarray}}
\newcommand{\eee}{\end{eqnarray}}
\newcommand{\eb}{\begin{equation}}
\newcommand{\ee}{\end{equation}}
\newcommand{\bT}{{\bf T}}
\newcommand{\bJ}{{\bf J}}
\newcommand{\barcalU}{{\bar{\cal{U}}}}
\newcommand{\barcalE}{{\bar{\cal{E}}}}
\newcommand{\pcalU}{{\bar{\cal{U}}}}
\newcommand{\pcalE}{{\bar{\cal{E}}}}
\newcommand{\calW}{{\cal{W}}}
\newcommand{\calP}{{\cal{P}}}
\newcommand{\calPma}{\calP_{\max}}
\newcommand{\calPmi}{\calP_{\min}}
\newcommand{\calPst}{\calP_{\rm sta}}
\newcommand{\calPi}{\calP_{\inf}}
\newcommand{\calPs}{\calP_{\sup}}
\newcommand{\calPvi}{\calP_{\rm vi}}
\newcommand{\calPe}{\calP_{\rm ext}}
\newcommand{\calPbv}{\calP_{\rm bv}}
\newcommand{\calPcbv}{\calP_{\rm cbv}}
\newcommand{\calPiv}{\calP_{\rm iv}}
\newcommand{\calLe}{L_{\rm ext}}
\newcommand{\calPm}{\calP_{\min}}
\newcommand{\Pbv}{\calP_{\rm bv}}
\newcommand{\calPbvc}{\calP_{\rm cbv}}
\newcommand{\Bvp}{\calP_{\rm vp}}

\newcommand{\sT}{\cal{T}}
\newcommand{\dsT}{\dot{\cal{T}}}
\newcommand{\bg}{{\bf g}}
\newcommand{\bm}{{\bf m}}
\newcommand{\bC}{{\bf C}}
\newcommand{\bG}{{\bf G}}
\newcommand{\bH}{{\bf H}}
\newcommand{\bn}{{\bf n}}
\newcommand{\bt}{{\bf t}}
\newcommand{\bff}{{\bf f}}
\newcommand{\bV}{{\bf V}}
\newcommand{\bY}{{\bf Y}}
\newcommand{\BR}{{\bf R}}
\newcommand{\bR}{{\bf R}}
\newcommand{\bv}{{\bf v}}
\newcommand{\bw}{{\bf w}}
\newcommand{\BN}{{\bf N}}
\newcommand{\bx}{{\bf x}}
\newcommand{\by}{{\bf y}}
\newcommand{\BX}{{\bf X}}
\newcommand{\bX}{{\bf X}}
\newcommand{\BF}{{\bf F}}
\newcommand{\bN}{{\bf N}}
\newcommand{\bF}{{\bf F}}
\newcommand{\bA}{{\bf A}}
\newcommand{\bB}{{\bf B}}
\newcommand{\bI}{{\bf I}}
\newcommand{\bj}{{\bf j}}
\newcommand{\bk}{{\bf k}}
\newcommand{\calS}{{\cal S}}
\newcommand{\calSs}{{\calS_s}}
\newcommand{\calSa}{{\calS_a}}
\newcommand{\calM}{{\cal M}}
\newcommand{\calMs}{{\calM_{sym}}}
\newcommand{\calMp}{\calM_+}
\newcommand{\calMo}{\calM_{ort}}
\newcommand{\DelM}{{\Delta_\calM}}
\newcommand{\DelMo}{{\Delta_\calM^o}}
\newcommand{\calbM}{\bar{\cal M}}
\newcommand{\calF}{{\cal F}}
\newcommand{\calG}{{\cal G}}
\newcommand{\calB}{{\cal B}}
\newcommand{\calb}{{\cal b}}
\newcommand{\calD}{{\cal D}}
\newcommand{\calH}{{\cal H}}
\newcommand{\calL}{{\cal L}}
\newcommand{\calK}{{\cal K}}
\newcommand{\calC}{{\cal C}}
\newcommand{\calT}{{\cal T}}
\newcommand{\calZ}{{\cal Z}}
\newcommand{\calU}{{\cal U}}
\newcommand{\calE}{{\cal E}}
\newcommand{\calEc}{{\calE_c}}
\newcommand{\calN}{{\cal N}}
\newcommand{\calV}{{\cal V}}
\newcommand{\calR}{{\cal R}}
\newcommand{\calY}{{\cal Y}}
\newcommand{\calX}{{\cal X}}
\newcommand{\dcalV}{\dot{\cal V}}
\newcommand{\dcalU}{\dot{\cal U}}
\newcommand{\dcalE}{{\dot{\cal E}}}
\newcommand{\dbq}{{\dot{\bf q}}}
\newcommand{\bu}{{\bf u}}
\newcommand{\bara}{\bar{ a}}
\newcommand{\barba}{\bar{\bf a}}
\newcommand{\barbE}{\bar{\bf E}}
\newcommand{\barbC}{\bar{\bf C}}
\newcommand{\barbF}{\bar{\bf F}}
\newcommand{\barbM}{\bar{\bf M}}
\newcommand{\barbN}{\bar{\bf N}}
\newcommand{\barbS}{\bar{\bf S}}
\newcommand{\barbT}{\bar{\bf T}}
\newcommand{\barbn}{\bar{\bf n}}
\newcommand{\barbp}{\bar{\bf p}}
\newcommand{\barbq}{\bar{\bf q}}
\newcommand{\barbu}{\bar{\bf u}}
\newcommand{\barf}{\bar{f}}
\newcommand{\barbf}{\bar{\bf f}}
\newcommand{\barbr}{\bar{\bf r}}
\newcommand{\barbx}{\bar{\bf x}}
\newcommand{\barby}{\bar{\bf y}}
\newcommand{\barm}{{\bar{m}}}
\newcommand{\barw}{{\bar{w}}}
\newcommand{\barb}{{\bar{b}}}
\newcommand{\bare}{{\bar{e}}}
\newcommand{\barimath}{{\bar{\imath}}}
\newcommand{\barbb}{{\bar{\bf b}}}
\newcommand{\barS}{{\bar{S}}}
\newcommand{\bart}{{\bar{t}}}
\newcommand{\barbt}{{\bar{\bf t}}}
\newcommand{\barbs}{{\bar{\bf s}}}
\newcommand{\barbv}{{\bar{\bf v}}}
\newcommand{\barnabla}{\bar{\nabla}}
\newcommand{\bgra}{{\bf \nabla}}
\newcommand{\dW}{{\dot{W}}}
\newcommand{\dU}{\dot{U}}
\newcommand{\dlam}{{\dot{\lambda}}}
\newcommand{\barlam}{{\bar{\lambda}}}
\newcommand{\bs}{{\bf s}}
\newcommand{\bb}{{{\bf b}}}
\newcommand{\bd}{{\bf d}}
\newcommand{\bl}{{\bf l}}
\newcommand{\barl}{{\bar{\bf l}}}
\newcommand{\dotbS}{\dot{\bf S}}
\newcommand{\dotbu}{\dot{\bf u}}
\newcommand{\dotbp}{\dot{\bf p}}
\newcommand{\dotbx}{\dot{\bf x}}
\newcommand{\dE}{\dot{\bf E}}
\newcommand{\bP}{{\bf P}}
\newcommand{\bS}{{\bf S}}
\newcommand{\dtau}{\,\mbox{d}\tau}
\newcommand{\dI}{{\,\mbox{d}I}}
\newcommand{\dt}{{\,\mbox{d}t}}
\newcommand{\ds}{{\,\mbox{d}s}}
\newcommand{\dd}{\mbox{d}}
\newcommand{\ddt}{\frac{\d}{\dt}}
\newcommand{\dbt}{{\dot{\bf t}}}
\newcommand{\dv}{\dot{\bf v}}
\newcommand{\dotv}{\dot{v}}
\newcommand{\du}{\dot{\bf u}}
\newcommand{\dotu}{\dot{ u}}
\newcommand{\dotxi}{\dot{\xi}}
\newcommand{\be}{{\bf e}}
\newcommand{\de}{\dot{\bf e}}
\newcommand{\dS}{\mbox{ d}S}
\newcommand{\dx}{\mbox{ d}x}
\newcommand{\ddx}{\frac{\d}{\dx}}
\newcommand{\dy}{\,\mbox{d}y}
\newcommand{\dO}{\,\mbox{d}\Oo}
\newcommand{\dOt}{\,\mbox{d}\Oot}
\newcommand{\dG}{\,\mbox{d} \Gamma}
\newcommand{\dbM}{{\rm d} \partial\calM}
\newcommand{\dM}{{\rm d} \calM}
\newcommand{\bU}{{\bf U}}
\newcommand{ \grad }{{\mbox{grad}}}
\newcommand{\alp}{{\alpha}}
\newcommand{\ab}{{{\alpha\beta}}}
\newcommand{ \eps}{{\epsilon}}
\newcommand{ \Ups}{{\Upsilon}}
\newcommand{ \sig}{{\sigma}}
\newcommand{ \Lam}{{\Lambda}}
\newcommand{ \barLam}{\bar{\Lambda}}
\newcommand{ \Gam}{{\Gamma}}
\newcommand{ \lam}{{\lambda}}
\newcommand{ \Xx}{{\bf X}}
\newcommand{ \Ff }{{\bf F }}
\newcommand{\barR }{{\bar{\real}}}
\newcommand{\barOo }{{\bar{\Omega}}}
\newcommand{\bc}{{\bf c}}
\newcommand{\bD}{{\bf D}}
\newcommand{\bZ}{{\bf Z}}
\newcommand{\epi}{{\mbox{epi }}}
\newcommand{\dom}{{\mbox{dom }}}
\newcommand{\Int}{{\mbox{int}}}
\newcommand{\tr}{{\mbox{tr}}}
\newcommand{\Lin}{{\cal{M}}}
\newcommand{\Cof}{{\mbox{Cof }}}
\newcommand{\rank}{{\mbox{rank }}}
\newcommand{\Diag}{{\mbox{Diag }}}
\newcommand{\Ker}{{\mbox{Ker }}}
\newcommand{\lin}{{\mbox{lin}}}

\newcommand{\Pis }{\Pi_{\sup}}
\newcommand{\Pist }{\Pi_{\small{\sta}}}
\newcommand{\Pii }{\Pi_{\inf}}

\newcommand{\com}{{\mbox{com }}}
\newcommand{\ine}{{\mbox{ine }}}
\newcommand{\Vol}{{\mbox{Vol}}}
\newcommand{\Inv}{{\mbox{Inv}}}
\newcommand{\dive}{{\mbox{div}}}
\newcommand{\curl}{{\mbox{curl}}}
\newcommand{\Sym}{{\mbox{Sym}}}
\newcommand{\Global}{{\mbox{Global}}}
\newcommand{\sym}{{\mbox{sym}}}
\newcommand{\divM}{{\mbox{div}}_\calM}
\newcommand{\VtoR}{\calU \rightarrow \real}
\newcommand{\mboxin}{${\mbox{ in }}$}
\newcommand{\mboxon}{{\mbox{ on }}}
\newtheorem{lemma}{Lemma}
\newtheorem{Cor}{Corollary}
\newtheorem{assumption}{Assumption}
\newtheorem{remark}{Remark}
\newtheorem{algorithm}{Algorithm}
\newtheorem{thm}{Theorem}
\newtheorem{Observation}{Observation}
\newtheorem{definition}{Definition}
\newtheorem{problem}{Problem}
\newtheorem{example}{Example}

   \newcommand\realmp{{\real^m_+}}
      \newcommand\realn{{\real_-}}
      \newcommand\realmn{{\real^m_-}}
      \newcommand\realp{{\real_+}}
      \renewcommand\pl{{|}}
        \renewcommand\calPs{{\calP_{\sta}}}
\renewcommand\AA{{A}}
\renewcommand\bA{{{A}}}
\renewcommand\bI{{{I}}}
\renewcommand\Pi{{{P}}}
\newcommand\BB{{B}}
\renewcommand\bb{{b}}
\renewcommand\bB{{B}}
\renewcommand\II{{\calI}}
\newcommand\FF{{F}}
\newcommand\yy{{y}}
\newcommand\byd{{y^*}}
\newcommand\KK{{K}}
\newcommand\y{{y}}
\renewcommand\by{{y}}
\newcommand\xxd{{{x^*}}}
\renewcommand\barbx{{\bar{x}}}
\renewcommand\bff{{{f}}}
\renewcommand\bc{{{f}}}
\newcommand\yyd{{y^*}}
\newcommand\baryd{{\bar{y}^*}}
\newcommand\barWd{{\barW^*}}
\newcommand\Wv{{\check{W}}}
\newcommand\calYd{\calY^*}
\newcommand\calXd{\calX^*}
\newcommand\G{{G\^{a}teaux} }

\renewcommand\eb{\begin{equation}}
\renewcommand\ee{\end{equation}}

\newcommand\calPd{{\calP^d}}
\newcommand\bepsd{{\beps^*}}
\newcommand\rhod{{\rho^*}}
\newcommand\PPd{{P^d}}
\newcommand\barbepsd{{\bar{\beps}^*}}
\newcommand\barxm{{{\bar{\bx}_\mu}}}
\newcommand\barhod{{\bar{\rho}^*}}
\renewcommand\calXk{{{\calX_f}}}
\newcommand\barbnu{{\bar{\bnu}}}
\renewcommand\AA{{A}}
\newcommand\QQ{{Q}}
\newcommand\XI{{\Xi}}
\newcommand\DD{{D}}
\newcommand\RR{{R}}
 \renewcommand\SS{{S}}
\newcommand\bq{\bf{q}}
\newcommand\barblam{\bar{\blam}}

\newcommand\GG{G}
\newcommand\GGm{G_{\mu}}
\newcommand\WW{{W}}
\renewcommand\bA{{{A}}}
\renewcommand\bI{{{I}}}
\renewcommand\Pi{{{P}}}
\renewcommand\bb{{\bf b}}
\renewcommand\bB{{B}}
\renewcommand\II{{\calI}}
\newcommand\xxo{{{\bf x}^o}}
\renewcommand\barbx{{\bar{\bf x}}}
\renewcommand\yy{{{\bf y}}}
\renewcommand\zz{{{\bf z}}}
\renewcommand\ww{{{\bf w}}}
\renewcommand\barww{{{\bar{\ww}}}}
\renewcommand\bc{{{\bf c}}}
\renewcommand\barbx{{\bar{\xx}}}
\newcommand\Nb{N_{A}}
\renewcommand\bff{{{\bf f}}}
\renewcommand\gg{g}

\title{\textbf{Complete Solutions to Nonconvex Fractional Programming Problems}}
\author{David Yang Gao $^{1}$\thanks{Email: d.gao@ballarat.edu.au}$\;\;$ and  Ning Ruan $^{1,2}$\thanks{Email: n.ruan@ballarat.edu.au}}
\date{}
\maketitle

\begin{center}
{\textbf{ David Yang Gao$^{1}$ and Ning Ruan$^{1,2}$}}\\
1. School of Sciences, Information Technology and Engineering,\\
University of Ballarat, Ballarat, Vic 3353, Australia.\\
2. Department of Mathematics and  Statistics, \\
Curtin University, Perth,  WA 6845, Australia.
\end{center}

\begin{abstract}
This paper presents a canonical dual approach to the problem of
minimizing the sum of a quadratic function and the ratio of nonconvex
function and quadratic functions, which is a type of non-convex optimization
problem subject to an elliptic constraint.
We first relax the fractional structure by introducing a family of
parametric subproblems. Under certain conditions,
we show that the canonical dual of each subproblem
becomes a two-dimensional concave maximization problem that exhibits
no duality gap. Since the infimum of the optima of the parameterized
subproblems leads to a solution to the original problem, we then
derive some optimality conditions and existence conditions for
finding a global minimizer of the original problem.

\end{abstract}

\pagebreak

\section{Introduction}
We study in this paper the following nonconvex fractional
programming problem:
\eb
(\calP ):\;\;\; \ \ \min \left\{ P_0(\bx)=  f(\bx) + \frac{g(\bx)}{h(\bx)}
\;\;  :  \;\;\; {\bx \in \calX} \right\}, \label{eq:govp}
\ee
where $\bx= ( x_1 , x_2, \cdots , x_n)^T \in \real^n$ and
\[
f(\bx)= \half \bx^T \QQ \bx  - \bff^T \bx,\;\;
\gg(\bx) = \half(\half|\BB \bx|^2-\lambda)^2 , \;\;
h(\bx) = \half \bx^T H \bx - \bb^T \bx,
\]
with $\BB \in \real^{m\times n }$, $\QQ  \in \real^{n\times n} $ being symmetric,
$H\in\real^{n\times n}$ negative definite , $\lam\in \real^+$, and $\bff, \bb \in
\real^n$, where $|v|$ denotes the Euclidean
norm of $v$. Assume that $\mu_0^{-1}=h(H^{-1}\bb)>0$ and
$\delta\in(0,\mu_0^{-1}]$, then the feasible domain $\calX$ is
defined by
\[
\calX = \{ \bx \in \real^n \; | \;\; h(\bx)  \ge \delta > 0 \},
\]
which is a constraint of elliptic type.

Problem $(\calP)$ belongs to a class of ``sum-of-ratios'' problems
that have been actively studied for several decades. The ratios
often stand for efficiency measures representing
performance-to-cost, profit-to-revenue or return-to-risk
for numerous applications in economics,
transportation science, finance, engineering, etc.
\cite{71Almogy,69Colantoni,92Falk,81Kanchan,89Konno,96Konno,95Schaible_handbook,Stancu-1980book}.
Depending on the nature of each application, the functions $f,g,h$
can be affine, convex, concave, or neither. However, even for the
simplest case in which $f,g,h$ are all affine functions, problem
$(\calP)$ is still a global optimization problem that may have
multiple local optima \cite{02Cambini,77Schiable}. In particular,
Freund and Jarre \cite{01FreundJarre} showed that the sum-of-ratios
problem $(\calP)$ is NP-complete when $f,g$ are convex and $h$ is concave.

Due to the non-convexity involved in the fractional structure, the
ordinary Lagrangean dual only provides a weak duality theorem that
may exist a positive duality gap. In this paper, we
explore some interesting properties and develop a canonical dual approach based on
Gao's work \cite{gao-book00} for solving problem $(\calP)$.

In Section 2, we first parameterize problem $(\calP)$ into a family
of subprograms $\{(\calP_\mu)\}$, in which each subproblem is a
non-convex quadratic program subject to one quadratic
constraint. Then, we show the infimum of the optima of the
parameterized subproblems provides a solution to problem $(\calP)$.
Since each subproblem $(\calP_\mu)$ is a non-convex problem, a
canonical dual problem $(\calP^d_\mu)$ is derived. We provide some
sufficient conditions to establish both the weak and strong duality
theorems (the so called {\it perfect duality}) for the pair of
$(\calP_\mu)$ and $(\calP^d_\mu)$. In Section 3, we develop some existence conditions
under which a global optimizer of the original problem $(\calP)$ can
indeed be identified by solving the corresponding canonical dual
problems.

\section{Sufficiency for Global Optimality}

In order to solve problem $(\calP)$, we consider the following
family of parameterized subproblem:
\eb
(\calP_\mu):\;\;\;   \min \left\{ \PP_{\mu}(\bx) = \half \bx^T
\QQ \bx  - \bff^T \bx  +  \mu  g(\bx) \;\;  :
\;\;\;{\bx \in \calX_\mu} \right\}, \label{problem_pmu}
\ee
where $\mu\in[\mu_0,\delta^{-1}]$ and
\[
\calX_\mu = \{ \bx \in \real^n \; | \;\; h(\bx) \ge \mu^{-1} \ge \delta > 0 \}
\]
is a convex set. We immediately have the following result:

\begin{lemma}\label{equivalent_form}
Problem $(\calP)$ is equivalent to  $(\calP_\mu)$ in the sense that
\eb
\inf_{\bx \in \calX} \PP_0(\bx) =\inf_{\mu\in[\mu_0,\delta^{-1}]}\inf_{\bx \in \calX_\mu}
\PP_{\mu}(\bx).
\ee
\end{lemma}

{\bf Proof}. It is easy to see that
\begin{eqnarray*}
&\,&\inf_{\bx \in \calX} \PP_0(\bx)\\
&=&\inf\limits_{\bx\in\calX} \left\{ \half\bx^T \QQ \bx-\bff^T \bx
+\frac{\gg(\bx)}{h(\bx)} \right\} \\
&=&\inf\limits_{\mu\in[\mu_0,\delta^{-1}]}\;\;\inf\limits_{h(\bx)=
\mu^{-1}}\left\{\half\bx^T \QQ \bx-\bff^T \bx+\frac{\gg(\bx)} {h(\bx)} \right\}\\
&=&\inf\limits_{\mu\in[\mu_0,\delta^{-1}]}\;\; \inf\limits_{h(\bx)=
\mu^{-1}}\left\{ \half\bx^T \QQ \bx-\bff^T \bx+\mu{\gg(\bx)}\right\}\\
&\ge&\inf\limits_{\mu\in[\mu_0,\delta^{-1}]} \;\; \inf\limits_{
\bx\in\calX_\mu}\left\{\half\bx^T \QQ \bx-\bff^T \bx+\mu{\gg(\bx)}\right\} \\
&=&\inf_{\mu\in[\mu_0,\delta^{-1}]}\;\; \inf_{\bx \in \calX_\mu}\PP_{\mu}(\bx).\\
\end{eqnarray*}

Conversely,
\begin{eqnarray*}
&\,&\inf\limits_{\mu\in[\mu_0,\delta^{-1}]}
\;\;\inf\limits_{\bx\in\calX_\mu}
\left\{ \half\bx^T \QQ \bx-\bff^T \bx+\mu{\gg(\bx)}\right\}\\
&=&\inf \limits_{\mu\in[\mu_0,\delta^{-1}]}
\;\;\inf\limits_{h(\bx)\ge\mu^{-1}}
\left\{ \half\bx^T \QQ \bx-\bff^T \bx+\mu{\gg(\bx)}\right\} \\
&\ge&\inf\limits_{\mu\in[\mu_0,\delta^{-1}]}
\;\;\inf\limits_{h(\bx)\ge \mu^{-1}} \left \{ \half\bx^T \QQ
\bx-\bff^T \bx+\frac{\gg(\bx)}{h(\bx)} \right\}\
(\hbox{since }\gg(\bx)>0)\\
&=&\inf_{\bx \in \calX} \PP_0(\bx).\\
\end{eqnarray*}
This completes the proof of the lemma.
\hfill $\Box$\\

Now, for any $\mu\in[\mu_0,\delta^{-1}]$, we define
\begin{eqnarray}
\GG_\mu(\vsig,\sig)&=&\QQ +\mu \BB^T \BB\vsig -\sigma \HH,
for\; \vsig\geq -\lambda,\;\sigma \geq 0,\label{gmu}\\
\calS_\mu^+&=&\{\vsig\geq -\lambda,\;\sigma \geq 0| \;\GG_\mu(\vsig,\sig) \succ 0\},
\end{eqnarray}
where `$\succ$' means positive definiteness of a matrix.
Then, the parametrical canonical dual problem can be proposed as the following:
\begin{eqnarray}
\PP_{\mu}^d(\vsig, \sigma)= -\half(\bff-\sigma \bb)^T
\GG^{-1}_{\mu} (\vsig, \sigma)(\bff -\sigma \bb)
-\mu \lambda \vsig-\frac{\mu}{2}\vsig^2+ \frac{\sig}{\mu}.
\end{eqnarray}

Given any $\mu\in[\mu_0,\delta^{-1}]$, consider the following
canonical dual problem $(\calP^d_\mu )$:
\[
\begin{array}
[c]{ccc}%
(\calP^d_\mu)\ \ \ \ &\sup & \PP^d_\mu(\vsig, \sig)  \\
&\rm {s.t.} & (\vsig, \sig)\in \calS^+_\mu.
\end{array}%
\]

\begin{thm}(Weak Duality)\label{weak_duality}
If there exists a global maximizer $(\vsig_{\mu},\sig_{\mu})$
of $P_{\mu}^d(\vsig, \sig)$ over $\calS_{\mu}^+$,
then the vector
\eb
\bx_{\mu}= G_{\mu}^{-1}(\vsig_{\mu,},\sig_{\mu})(\bff-\sig_{\mu} \bb)
\ee
is a global minimizer of $(\calP_{\mu})$ over $\calX_{\mu}$ and
\eb
P_{\mu}^d(\vsig,\sig)\leq P_{\mu}(\bx),\;
\forall (\bx,\vsig,\sig) \in \calX_{\mu}\times \calS_{\mu}^+.
\ee
\end{thm}
{\bf Proof}.
Let $\Lam (\cdot ): \real^n \rightarrow \real $
be the \textit{geometrical transformation}\cite{gao-book00,gao-sherali07,gao-strang89} defined by
\eb
\xi= \Lam (\bx) =
\half |\BB \bx |^2 - \lambda
\ee
and let
\eb
\UU(\xi)=\half \veps^2
\ee
Then, Problem $(\calP_\mu)$ in (\ref{problem_pmu}) can be written as
the following unconstrained optimization problem
\eb\label{uncon}
\min\left\{P(\bx)=\half \bx^T \QQ \bx -\bff^T \bx
+\mu(\UU(\Lambda(\bx))- \sigma(h(\bx)-\mu^{-1})
|\bx \in \real^n.\right\}.
\ee
Let $\vsig$ be the dual variable of $\xi$, i.e.,
$\vsig=\nabla \UU(\xi)=\xi$,
the Legendre conjugate $\UU^*(\vsig)$ can be uniquely defined by
\eb
\UU^\sharp(\vsig)=\sta_{\xi \geq \lam}\{\xi \vsig- \UU(\xi)\}= \half \vsig^2
\ee
where $\vsig \in \calV_a^*= \{\vsig\in \real|\vsig\geq -\lambda\}.$

By replacing $\UU(\Lambda(\bx))$ with
$ \Lambda(\bx)^T \vsig -\UU^{\sharp}(\vsig)$ in (\ref{uncon}),
we define the \textit{total complementary function} as
\begin{eqnarray*}
\Xi(\bx,\vsig,\sig) &=& \half \bx^T \QQ\bx-\bff^T \bx+ \mu (\UU(\Lambda(\bx))
-\sig(\half \bx^T \HH \bx- \bb^T \bx- \mu^{-1})\\
&=& \half \bx^T \QQ \bx -\bff^T \bx+ \mu( \Lambda
(\bx)^T \vsig- \UU^{\sharp} (\vsig))
-\sig(\half \bx^T \HH\bx -\bb^T \bx-\mu^{-1})\\
&=& \half \bx^T \QQ \bx -\bff^T \bx+ \mu[(\half(|\BB\bx|^2-\lambda)\vsig-\half \vsig^2]
-\sig(\half \bx^T\HH \bx-\bb^T \bx -\mu^{-1})\\
&=& \half \bx^T \GG_{\mu}(\vsig,\sig)\bx -(\bff-\sig \bb)^T \bx
-\mu \lambda \vsig-\frac{\mu}{2}\vsig^2 +\frac{\sig}{\mu},
\end{eqnarray*}
where $G_\mu(\vsig,\sig)$ is defined in (\ref{gmu}). Note that $\XI(\bx,\vsig,
\sig)$ is convex in $\bx\in\real^n$ for any given
$(\vsig,\sig)\in\calS^+_\mu$ and affine (hence concave) in $(\vsig,\sig)$ for any
given $\bx\in\real^n$. By the criticality condition
\eb
\frac{\partial \XI}{\partial{\bx}}=G_\mu(\vsig, \sig)\bx-(\bff-\sig\bb)= 0.\label{primal_critical}
\ee
we have $\bx(\sig)=\GGm^{-1}(\vsig, \sig)(\bff-\sig \bb)$, which is the
global minimizer of $\XI(\bx,\sigma)$. Moreover,
\begin{eqnarray*}
\min_{\bx\in\real^n}\XI(\bx,\vsig,\sig)&=&\XI(\bx(\vsig,\sig),\vsig,\sig)\\
&=&\half\bx(\vsig,\sig)^T(G_\mu(\vsig,\sig))\bx(\vsig,\sig)-
(\bff-\sig\bb)^T\bx(\vsig,\sig)-\mu \lambda \vsig-\frac{\mu}{2}\vsig^2 +\frac{\sig}{\mu}\\
&=&\half\bx(\vsig,\sig)^T(\bff-\sig\bb)
-(\bff-\sig\bb)^T\bx(\vsig,\sig)-\mu \lambda \vsig-\frac{\mu}{2}\vsig^2 +\frac{\sig}{\mu}\\
&=&-\half(\bff-\sig\bb)^T\bx(\vsig,\sig)-\mu \lambda \vsig-\frac{\mu}{2}\vsig^2 +\frac{\sig}{\mu}\\
&=&-\half(\bff-\sig\bb)^T
G_\mu^{-1}(\vsig,\sig)(\bff-\sig\bb)-\mu \lambda \vsig-\frac{\mu}{2}\vsig^2 +\frac{\sig}{\mu}\\
&=&\PP^d_\mu(\vsig,\sig).
\label{saddle_1}
\end{eqnarray*}

By the assumption, $(\vsig_{\mu},\sig_{mu})$ is a global maximizer of $P_{\mu}^d(\vsig,\sig)$.
If $(\vsig_{\mu},\sig_{\mu})$ is an interior of $\calS_{\mu}^+$, then
$\frac{\partial}{\partial \vsig}P_{\mu}^d(\vsig_{\mu},\sig_{\mu})=0$, and
$\frac{\partial}{\partial \sig}P_{\mu}^d(\vsig_{\mu},\sig_{\mu})=0$. Otherwise,
we have $\frac{\partial}{\partial \vsig}P_{\mu}^d(\vsig_{\mu},\sig_{\mu})=0$,
$\sig_{\mu}=0$, $\frac{\partial}{\partial \sig}P_{\mu}^d(\vsig_{\mu},\sig_{\mu})\leq 0$.
In either case, If we denote $\bx_{\mu}=\bx(\vsig_{\mu},\sig_{\mu})=
G_{\mu}^{-1}(\vsig_{\mu},\sig_{\mu})(\bff-\sig_{\mu} \bb)$, we have
\begin{eqnarray*}
\frac{\partial}{\partial \vsig}P_{\mu}^d(\vsig_{\mu},\sig_{\mu})&=&
\mu(\half \bx(\vsig_{\mu},\sig_{\mu})^T B^T B \bx(\vsig, \sig)-\lam -\vsig)=0,\\
\frac{\partial}{\partial \sig}P_{\mu}^d(\vsig_{\mu},\sig_{\mu})&=&
\frac{1}{\mu}-\bx_{\mu}^T (\half H \bx_{\mu}-\bb)\leq 0.
\end{eqnarray*}
That is,
\begin{eqnarray*}
&&\vsig = \half | B\bx|^2-\lam,\\
&&\half \bx^T H \bx -\bb^T \bx -\mu^{-1} \geq 0.
\end{eqnarray*}
Therefore, $\bx_{\mu} \in \calX_{\mu}$, and for any $(\vsig,\sig) \in \calS_{\mu}^+$,
we have
\begin{eqnarray*}
P_{\mu}^d(\vsig,\sig) &\leq & P_{\mu}^d(\vsig_{\mu}, \sig_{\mu})\\
&=& \min_{\bx \in \real^n} \Xi( \bx, \vsig_{\mu},\sig_{\mu})\\
&=& \Xi(\bx_{\mu},\vsig_{\mu},\sig_{\mu})\\
&=& \min_{\bx \in \calX_{\mu}} \Xi(\bx, \vsig_{\mu},\sig_{\mu})\\
&=& \half \bx^T \QQ \bx-\bff^T \bx
+\mu(\Lam(\bx)^T \vsig-U^{\sharp} (\vsig))
-\sig(\half \bx^T H \bx-\bb^T \bx-\mu^{-1})\\
&\leq& \half \bx^T \QQ \bx-\bff^T \bx
+\mu(\Lam(\bx)^T \vsig-U^{\sharp} (\vsig))\\
&=& \half \bx^T \QQ \bx-\bff^T \bx
+ \mu(\half(\half|Bx|^2 -\lam)^2=P_{\mu}(\bx).
\end{eqnarray*}
This completes the proof.
\hfill $\Box$

\begin{thm}(Strong Duality)\label{strong_duality}
If $(\vsig_{\mu},\sig_\mu)$ is a critical point of $\PP^d_\mu(\vsig,\sig)$ over
$\calS^+_\mu$, then $(\calP^d_\mu)$ is perfectly dual to
$(\calP_\mu)$ in the sense that the vector
\eb
\bx_\mu =\GGm^{-1}(\vsig_{\mu}, \sig_\mu)(\bff-\sig_\mu \bb)
\ee
is a global minimizer of
$(\calP_\mu)$ and $(\vsig_{\mu}, \sig_\mu) $ is a global maximizer
of $(\calP^d_\mu)$, and
\begin{eqnarray}
\min_{\bx \in \calX_\mu}\PP_\mu(\bx) = \PP_\mu(\bx_\mu )
= \PP^d_\mu(\vsig_{\mu},\sig_\mu )
= \max_{(\vsig,\sig) \in \calS^+_\mu} \PP^d_\mu(\vsig,\sig).
\label{eq:minmax}
\end{eqnarray}
\end{thm}
{\bf Proof}.
The proof basically follows that of the former weak duality Theorem, The only
difference lies in the assumption that $(\vsig_{\mu},\sig_{\mu})$ is
a critical point of $P_{\mu}^d(\vsig,\sig) $ over $\calS_{\mu}^+$. In this case,
$\frac{\partial}{\partial \vsig}P_{\mu}^d(\vsig_{\mu},\sig_{\mu})=0$, and
$\frac{\partial}{\partial \sig}P_{\mu}^d(\vsig_{\mu},\sig_{\mu})=0$.
So $\bx_{\mu}=\bx(\vsig_{\mu},\sig_{\mu})=
G_{\mu}^{-1}(\vsig_{\mu},\sig_{\mu})(\bff-\sig_{\mu}\bb)$ is on the boundary of $\calX_{\mu}$.
That is, $\half \bx^T H \bx -\bb^T \bx-\mu_{-1}=0$. this further implies that
\eb
P_{\mu}^d(\vsig_{\mu},\sig_{\mu})=
\Xi(\bx_{\mu},\vsig_{\mu},\sig_{\mu})=
P_{\mu}(\bx_{\mu})
\ee
and the equation (\ref{eq:minmax}) follows naturally.
\hfill $\Box$

The above results immediately lead to the following sufficient
condition for finding the global optimizer of problem $(\calP)$:
\begin{Cor}
If $(\vsig_\mu,\sig_\mu)\in\calS^+_\mu$ holds for all
$\mu\in[\mu_0,\delta^{-1}]$, then
\eb
\min_{\bx\in\calX}
P_0(\bx)=\min_{\mu\in[\mu_0,\delta^{-1}]} \PP^d_\mu(\vsig_\mu,\sig_\mu ).
\ee
\end{Cor}

\section{Existence of Global Optimality }

Before we provide the condition for the existence of a global optimal solution
$(\vsig_{\mu},\sig_{\mu})$ to problem $(\PP_{\mu}^d)$ over $\calS_{\mu}^+$ with
any given $\mu\in[\mu_0,\delta^{-1}]$,
we need the following property.
\begin{lemma}\label{monotone}
For any $\mu\in[\mu_0,\delta^{-1}]$,
$\PP^d_{\mu}(\sig, \vsig)$ is a two-dimensional concave function over $\calS^+_\mu$.
\end{lemma}
{\bf Proof}.
Notice that the Hessian Matrix of the dual objective function is
\begin{eqnarray*}
\nabla^2 \PPd (\sig,\vsig)=S
= \left( \begin{array}{cc}
H_{\sig^2}& H_{\sig\vsig}\\
H_{\vsig\sig}& H_{\vsig^2} \end{array}
\right),
\end{eqnarray*}
where
\begin{eqnarray*}
H_{\sig^2}&=& -(H x(\vsig,\sig)-b) G^{-1}(\bsig,\sig)(H\bx(\bsig,\sig)-b),\\
H_{\vsig^2}&=& -(\mu B^T B)^2 G^{-1}(\vsig,\sig)\bx(\vsig,\sig)-\mu I,\\
H_{\vsig\sig}&=& \mu B^T B \bx(\vsig,\sig) G^{-1}(\vsig,\sig)(H \bx(\vsig,\sig)-b).
\end{eqnarray*}
In order to show the dual function is a concave function, it is equivalent to show that
\begin{eqnarray*}
S_0
= \left( \begin{array}{cc}
H_{\sig^2}& H_{\sig\vsig}\\
H_{\vsig\sig}& H_{\vsig^2}+\mu \end{array}
\right),
\end{eqnarray*}
is semi-negative definite.
By Sylvester's Criterion, it suffices to show that all the leading principal minors have a
non-positive determinant. Obviously, the first $n-1$ leading principal minors have
non-positive determinants, since
\eb
-(H \bx(\vsig,\sig)-b) G^{-1}(\vsig,\sig)(H\bx(\vsig,\sig)-b)
\ee
is semi-negative definite. It is left to show $det(S_0)\leq 0$. Note that
\begin{eqnarray*}
S_0 &=& C\cdot D\\
&=& \left( \begin{array}{cc}
-(H\bx(\vsig,\sig)-b)G^{-1}(\vsig,\sig)& 0\\
0& U B^T BG^{-1}(\vsig,\sig)\bx(\vsig,\sig) \end{array}
\right)\\
&&\cdot
\left( \begin{array}{cc}
H \bx(\vsig,\sig)-b& U B^T B\\
H \bx(\vsig,\sig)-b& U B^T B\end{array}
\right)
\end{eqnarray*}
Apparently, $Rank(C \cdot D) \leq Rank (D)\leq n $. We can make a conclusion that
$det S_0=0$. Thus, $S$ is semi-negative definite, which implies that dual function
is concave function.
\hfill $\Box$\\

Let $\partial \calS_{\mu}^+$ denotes a singular hyper-surface defined by
\eb
\partial \calS_{\mu}^+= \{ \bsig \geq -\lam, \sig \geq 0| \;\;
G_{\mu}(\vsig,\sig)  \succeq 0 , \;\; \det G_{\mu}(\vsig,\sig) =  0 \}.
\ee
\begin{thm}(Existence)\label{existence}
Given any $\mu\in[\mu_0,\delta^{-1}]$, if
\eb
\lim_{ (\vsig,\;\sig) \rightarrow \partial \calS_a^+}
\PP^d_{\mu}( \vsig,\sig ) = -\infty, \;\;  \forall (\vsig,\sig) \in \calS_{\mu}^+,
\label{condition1}
\ee
and
\eb
\lim_{ \vsig \rightarrow \infty,\; \sig \rightarrow -\infty}
\PP^d_{\mu}( \vsig,\sig ) = -\infty, \;\;  \forall (\vsig,\sig) \in \calS_{\mu}^+,
\label{condition2}
\ee
then the canonical dual problem $(\calP^d_\mu )$ has at least
one global optimal solution $(\vsig_{\mu},\sig_{\mu})\in\calS^+_\mu$.
\end{thm}
{\bf Proof}.
It follows from (\ref{condition1}) and
(\ref{condition2}) that there exists one $(\vsig_{\mu},\sig_{\mu})\in\calS^+_\mu$
to be a critical point of $\PP^d_\mu(\vsig,\sig)$. By Theorem
\ref{strong_duality}, we know $(\vsig_{\mu},\sig_{\mu})$ is a global maximizer
of $(\calP^d_\mu )$. \hfill $\Box$

\section{Conclusions}

In this paper, we study a kind of problems with sum of a quadratic function
and the ratio of nonconvex function and quadratic function as its objective function. We
first parameterize such a problem into a family of subproblems. Then we
develop a corresponding canonical duality theory, both in weak and strong
duality form,  to handle each subproblem. Based on the properties of the subproblems,
we provide non only the extremality conditions for global optimality of
the original problem, but also existence conditions
to assure that the global optimal solutions of the primal problems
can indeed be found by solving a sequence of concave maximization problems.

\noindent{\bf Acknowledgement}:
This paper was partially supported by a grant (AFOSR FA9550-10-1-0487)
from the US Air Force Office of Scientific Research. Dr. Ning Ruan was
supported by a funding from the Australian Government
under the Collaborative Research Networks (CRN) program.

\end{document}